\documentclass[final, 12pt, 3p,times]{elsarticle}

\usepackage{amsmath,amssymb,amsfonts,bbm,textcomp,amsthm,mathrsfs}

\usepackage{varioref}

\usepackage{mathtools}
\usepackage{comment}

\usepackage[h]{esvect}

\usepackage{array}

\usepackage{listings}
\usepackage[activate]{pdfcprot}
\usepackage{pdfpages}
\usepackage{booktabs}

\usepackage{textcomp}
\usepackage{tabulary}

\usepackage{multicol}
\usepackage{stmaryrd}

\usepackage{bbm}
\usepackage{bigints}
\usepackage{lmodern}
\usepackage{xcolor}

\DeclareFontFamily{U}{mathx}{\hyphenchar\font45}
\DeclareFontShape{U}{mathx}{m}{n}{
      <5> <6> <7> <8> <9> <10>
      <10.95> <12> <14.4> <17.28> <20.74> <24.88>
      mathx10
      }{}
\DeclareSymbolFont{mathx}{U}{mathx}{m}{n}
\DeclareFontSubstitution{U}{mathx}{m}{n}
\DeclareMathAccent{\widecheck}{0}{mathx}{"71}
\DeclareMathAccent{\wideparen}{0}{mathx}{"75}

\newcommand{\im}{\ensuremath{\mathrm{i}}}

\newcommand{\C}{\mathbb{C}}
\newcommand{\E}{{\mathbb E}}

\newcommand{\N}{\mathbb{N}}
\newcommand{\Pa}{{\mathbb P}}

\newcommand{\R}{\mathbb{R}}

\newcommand{\Bcal}{{\mathcal B}}

\newcommand{\Fcal}{{\mathcal F}}
\newcommand{\Gcal}{{\mathcal G}}

\newcommand{\Jcal}{{\mathcal J}}

\newcommand{\Ucal}{{\mathcal U}}

\newcommand{\Ncal}{{\mathcal N}}

\newcommand{\Rcal}{{\mathcal R}}
\newcommand{\Scal}{{\mathcal S}}
\newcommand{\Tcal}{{\mathcal T}}

\newcommand{\Zcal}{{\mathcal Z}}

\newcommand{\Rp}[1]{{\mathbb R^{#1}_{\geq0}}}

\newcommand{\expvBig}[2]{\operatorname{\E}^{#1}\Bigl[#2\Bigr]}

\newcommand{\marry}{{\textrm{\textmarried}}}

\newcommand{\til}[2]{{\overset{\sim}{#1}}^{{\raisebox{-5pt}{\scriptsize{$#2$}}}}}
\newcommand{\notil}[2]{{\overset{\nsim}{#1}}^{{\raisebox{-5pt}{\scriptsize{$#2$}}}}}

\newcommand{\supp}{\mathop{\mathrm{supp}}}

\newcommand{\la}{\left\langle}
\newcommand{\ra}{\right\rangle}

\newcommand{\st}{\;|\;}

\renewcommand{\epsilon}{\ensuremath\varepsilon}

\renewcommand{\phi}{\ensuremath{\varphi}}

\SetSymbolFont{stmry}{bold}{U}{stmry}{m}{n}

\newtheorem{theorem}{Theorem}[section]
\newtheorem{example}[theorem]{Example}
\newtheorem{remark}[theorem]{Remark}

\newtheorem{definition}[theorem]{Definition}
\newtheorem{lemma}[theorem]{Lemma}
\newtheorem{proposition}[theorem]{Proposition}
\newtheorem{notation}[theorem]{Notation}

\newtheorem{assumption}[theorem]{Assumptions}

\definecolor{newtxt}{rgb}{1.0,0.01,0.24}

\journal{}

\begin{document}

\begin{frontmatter}



\title{Pathwise construction of affine processes}


\author{Nicoletta Gabrielli\fnref{nico}}
\author{Josef Teichmann \fnref{josef}}
\fntext[nico]{University of Z\"urich, Plattenstrasse 22, CH-8032, Switzerland. \texttt{nicoletta.gabrielli@bf.uzh.ch}}
\fntext[josef]{ETH Z\"urich, R\"amistrasse 101, CH-8092, Switzerland. \texttt{josef.teichmann@math.ethz.ch}}

\begin{abstract}
\noindent
Based on the theory of multivariate time changes for Markov processes, we show how to identify affine processes as solutions of certain time change equations. The result is a strong version of a theorem presented by J. Kallsen in \cite{didactic_2006} which provides a representation in law of an affine process as a time--change transformation of a family of independent L\'evy processes.

\end{abstract}



\end{frontmatter}


\setlength{\parindent}{0pt}

\section{Introduction}

During the last decades, many alternatives to the Black-Scholes model have been proposed in the literature to overcome its
deficiencies. Possible extensions include jumps, stochastic volatility and/or other high dimensional models. Among the most popular ones,
we recall the exponential L\'evy models, which generalize the Black-Scholes model by introducing jumps. They allow to generate
implied volatility smiles and skews similar to the ones observed in the markets. However, in some occasions, independence of
increments is too big a restriction. Stochastic volatility models give a way to overcome this problem: when we model the variance parameter in the Black--Scholes model by a CIR model, we get the Heston's model, see \cite{heston_closed-form_1993}. The Heston model can be extended by adding jumps in the return component, as in the Bates model (see \cite{bates_jumps_1996}), and also in the stochastic variance component, as in the Barndorff--Nielsen and Shephard model (see \cite{barndorff-nielsen_shepard}). The class of affine processes includes all the above mentioned examples.

Affine processes are a class of time homogeneous Markov processes $X=(X_t)_{t\geq0}$
taking values in a state space $D\subset\R^d$ characterized by the fact that, for all
$(t,x)\in\Rp{}\times D$, their characteristic
function has the following exponential affine form
$$\expvBig{x}{e^{\la u,X_t\ra}}=e^{\phi(t,u)+\la x,\psi(t,u)\ra},\qquad u\in\im\R^d\,,$$
where $\phi$ and $\psi$ are two function taking values in $\C$ and $\C^d$, respectively.
The theory of affine processes is dominated by weak characterizations, since affine processes are characterized by a property of their marginal distributions. The functions $\phi$ and $\psi$ in the specification of the affine property,
solve a system of ODEs, also known in the literature with the name of
\emph{generalized Riccati equations}. These equations arise from the regularity property of affine processes. More
precisely, in \cite{cuchiero_path_2011} it has been proved that, even on a general state space, stochastically continuous
processes having the aforementioned affine property admit a version with c\`adl\`ag trajectories. The path regularity implies that the process is a semimartingale with differentiable characteristics up to its lifetime. From this characterization it is possible to conclude differentiability with respect to time of the Fourier--Laplace transform. This property, also called regularity property, is crucial to relate the marginal laws of affine processes with a solution of a system generalized Riccati equations.

\bigskip

This paper is devoted to a pathwise construction of affine processes, when the state space is specified by $\Rp{m}\times\R^n$. The representation proposed in this paper is a multivariate generalization of the Lamperti transformation of L\'evy processes in $\R$ with no negative jump. When $D=\Rp{}$, it has been proved that there exists one-to-one correspondence between affine processes
taking values in $D$ and L\'evy processes, see \cite{caballero_lamperti-type_2013}. More precisely, let $Z^{(1)}=(Z^{(1)}_t)_{t\geq0}$ be a L\'evy process starting from $0$ taking values in $\R$,
whose L\'evy measure has support $\Rp{}$ and let $Z^{(0)}$ be an independent subordinator.
Theorem 2 in  \cite{caballero_lamperti-type_2013} shows that there exists a solution of the following time--change equation
$$X_t=x+Z^{(0)}_t+Z^{(1)}_{\int_0^t X_s ds}$$
for all $(t,x)\in\Rp{}\times \Rp{}$. Moreover, it is proved that the solution is a time homogeneous Markov process,
taking values in $\Rp{}$ starting from $x$, characterized by the property that the logarithm of the
characteristic function of the transition semigroup is given by an affine function of the initial state $x$.
Hence, by definition, it is an affine process taking values in $\Rp{}$.

In this paper we aim to obtain the analogous result in the multivariate case. In \cite{didactic_2006} it has been proved that -- in distribution -- affine processes can be represented by means of $d+1$ independent L\'evy processes taking values in $\R^d$. Under some natural assumptions on the L\'evy triplets, the time change equation
\begin{equation}\label{eq:multilam}
X_t=x+Z^{(0)}_t+\sum_{i=1}^d Z^{(i)}_{\int_0^t X^{(i)}_s ds},\qquad t\geq0\,,
\end{equation}
admits a \emph{weak solution}. More precisely, \eqref{eq:multilam} admits a weak solution if there exists a probability space $(\Omega,\Gcal, P)$ containing two processes
$( X, Z)$ such that \eqref{eq:multilam} holds.
\begin{remark}
Recall that, a priori, the process $X$ takes values on the state space $\Rp{m}\times\R^n$ and, therefore, for all $j=m+1,\ldots,d$, the process $X^{(j)}$ is real valued. As we announced, existence of a solution for \eqref{eq:multilam} will be proved under a set of conditions on the L\'evy triplets. In particular, we will require that, for all $j=m+1,\ldots,d$, the L\'evy process $Z^{(j)}$ is deterministic (see Table \ref{tab:admpar} at pag. \pageref{tab:admpar}). This ensures that the sum in \eqref{eq:multilam} is well defined also for the indices $i=m+1,\ldots,d$.
\end{remark}

However in \cite{didactic_2006} the following problem
is left unsolved: is $X$ a strong solution of the time change equation \eqref{eq:multilam}? In this paper we try to address this problem and the exposition
is organized as follows. In Chapter \ref{sec:pre} we provide an overview of
some basic results for affine processes. We introduce a particular class of affine processes, called of Heston type,
which, up to a pathspace transformation, represents the full class of affine processes. See Definition
\ref{def:hestonaff} and Proposition \ref{prop:gocanon}. Chapter \ref{sec:existstc} contains the core of the proof of existence of a strong solution
of \eqref{eq:multilam}. The final result is stated in Theorem \ref{teo:multilam} and the proof is divided in several steps. Using the results from Chapter \ref{sec:existstc}, we will see how to construct a solution $X$ of \eqref{eq:multilam} which lives on the same
probability space where the L\'evy processes are defined. In Chapter \ref{sec:pathconstr} we show that, starting from a family of L\'evy processes
$\{Z^{(k)}\}_{k=0,1,\ldots,d}$ specified by some restrictions on their L\'evy triplets, the solution
time--change equation \eqref{eq:multilam} is a time homogeneous Markov process having the affine property. Observe that, this new existence proof of affine processes gives, as straightforward consequence, the c\`adl\`ag
property for affine processes.

\section{Preliminaries}\label{sec:pre}
\subsection{Notation}

Henceforth $D$ denotes the subset $\Rp{m}\times\R^n$ of $\R^d.$ The canonical basis of $\R^d$ is denoted by $\{e_i\}_{i=1,\ldots,d}$.
Given $\Delta\notin D$ define $D_\Delta=D\cup\{\Delta\}$. The set $\Bcal(D)$ is the space of measurable function on $D$,
while $mbbd(D)$ is the space of measurable bounded function on $D$.

 In order to simplify the notation, we introduce the sets of indices $I$ and $J$
defined as
$$I=\{1,\ldots,m\}\qquad\mbox{ and }\qquad J=\{m+1,\ldots,d\}.$$
Moreover, given a set $H\subseteq\{1,\ldots,d\},$ the map $\pi_H$ is the projection of $\R^d$ on the lower dimensional subspace with indices in $H.$
In particular
\begin{eqnarray*}
{\pi_I}:\Rp{m}\times\R^n&\to&\Rp{m}\\
x&\mapsto&\pi_I x:=(x_i)_{i\in I}\\
\end{eqnarray*}
and
\begin{eqnarray*}
{\pi_J}:\Rp{m}\times\R^n&\to&\R^n\\
x&\mapsto&\pi_J x:=(x_j)_{j\in J}\, .\\
\end{eqnarray*}

Due to the geometry of the state space, the function
\begin{equation}\label{eq:fux}
f_u(x):=e^{\la x,u\ra}\,,\qquad x\in D
\end{equation}
is bounded if and only if
\begin{equation}\label{eq:ucal}
\Ucal:=\C^m_{\leq0}\times \im\R^n\, .
\end{equation}

The notation $\la\cdot,\cdot\ra$ with input variables in $\R^d$ denotes the usual scalar product. The same notation is used
also when the scalar product is considered in the space $\R^d+i \R^d$. In this case we mean the extension of $\la\cdot,\cdot\ra$
in $\R^d+i\R^d$ without conjugation.

Unless differently specified, the notation $\E^x[\cdot]$ indicates that the expectation is taken under the probability measure $\Pa^x\, .$

Fix $N\in\N$ and let $s\in\Rp{N}$. Whenever we are going to consider $s$ as a time parameter, we emphasize its multidimensionality by writing
$\underline s$. When $\underline s=(s_1,\ldots,s_N)$ is a multivariate time parameter and $X$ is a stochastic process in $\R^N$, we use the notation
$$\underline X(\underline s):=(X^{(1)}_{s_1},\ldots,X^{(N)}_{s_N})\in \R^{N} \, \, .$$

\subsection{Affine processes}

In line with the literature, we introduce the affine processes as a class of time homogeneous Markov processes characterized by two additional properties.
The first one being stochastic continuity, the second one a condition which characterizes the Fourier--Laplace transform of the one time marginal
distributions. This introduction of affine processes in taken from \cite{2003,cuchiero_path_2011} and \cite{keller-ressel_regularity_2011}.

\begin{definition}\label{def:affine} Let $$(\Omega, (X_t)_{t\geq0},(\Fcal^\natural_t)_{t\geq0},(p_t)_{t\geq0},(\Pa^x)_{x\in D})$$
be a time homogeneous Markov process. In particular we assume that
\begin{itemize}[-]
\item  $\Omega$ is a probability space,
\item $(X_t)_{t\geq0}$ is a stochastic process taking values in $D_\Delta,$
\item $\Fcal^\natural_t=\sigma(\{X_s\;, s\leq t\})$,
\item $(p_t)_{t\geq0}$ is a semigroup of transition functions on $(D_\Delta,\Bcal(D_\Delta)),$
\item $(\Pa^x)_{x\in D_\Delta}$ is a probability measures on $(\Omega,\Fcal^\natural)$,  with $\Fcal^\natural=\bigvee_{t\geq0}\Fcal^\natural_t$,
\end{itemize}
satisfying
\begin{equation}\label{eq:markovprop}
\expvBig{x}{f(X_{t+s})\big|\Fcal^\natural_t}=\expvBig{X_t}{f(X_s)},\quad \Pa^x\mbox{-a.s. for all }f\in mbdd(D_\Delta).
\end{equation}

The process $X$ is said to be an \emph{affine process} if it satisfies the following properties:
\begin{itemize}
 \item\label{prop:stoc}for every $t\geq0$ and $x\in D,$ $\lim_{s\to t}p_s(x,\cdot)=p_t(x,\cdot)$ weakly,
\item\label{prop:affine}there exist functions $\phi:\Rp{}\times\Ucal\to\C$ and $\psi:\Rp{}\times\Ucal\to \C^d$ such that
\begin{equation}\label{eq:affineprop}
\expvBig{x}{e^{\la u,X_t\ra}}=\int_D e^{\la u,\xi\ra}p_t(x,d\xi)=e^{\phi(t,u)+\la x,\psi(t,u)\ra},
\end{equation}
 for all $x\in D$ and $(t,u)\in\Rp{}\times\Ucal.$
\end{itemize}
\end{definition}

Regularity is a key feature for an affine process.
It gives differentiability of the Fourier--Laplace transform with respect to time.
\begin{definition}
 An affine process $X$ is called \emph{regular} if, for every $u\in\Ucal$, the
derivatives
\begin{equation}\label{regular}
F(u):=\partial_t \phi(t,u)\Big|_{t=0},\qquad R(u):=\partial_t\psi(t,u)\Big|_{t=0}\,,
\end{equation}
exist for all $u\in\Ucal$ and are continuous in $$\Ucal_m=\left\{ u\in \C^d\st \sup_{x\in D}\Rcal e (\la u,x\ra)\leq m\right\}\,,$$
for all $m\geq1$.
\end{definition}

Regularity has been proved in \cite{cuchiero_path_2011} for the class of affine processes on general state spaces. The proof is based on the fact that affine processes always admit a version which has c\`adl\`ag
paths. From this path regularity it is possible to conclude differentiability of the Fourier--Laplace transform. We summarize here the main results.

\begin{theorem}[Theorem 6.4 in \cite{cuchiero_path_2011}]\label{teo:affine} Every affine process is regular. On the set $\Rp{}\times\Ucal,$
the functions $\phi$ and $\psi$ satisfy the following system of generalized Riccati equations:
\begin{equation}\label{eq:riccati}
\begin{aligned}
   \partial_t\phi(t,u)&=F(\psi(t,u)),&\phi(0,u)=0\,,\\
\partial_t\psi(t,u)&=R(\psi(t,u)),&\psi(0,u)=u,
\end{aligned}
\end{equation}
with
\begin{eqnarray}
F(u)&=&\la b,u\ra+\frac{1}{2}\la u,au \ra-c\nonumber\\
&&+\int_{D\setminus\{0\}}\left(e^{\la u,\xi\ra}-1-\la \pi_J u,\pi_J h(\xi)\ra\right)m(d\xi),\label{eq:F}\\
R_k(u)&=&\la \beta_k,u\ra+\frac{1}{2}\la u,\alpha_ku\ra-\gamma_k\nonumber\\
&&+\int_{D\setminus\{0\}}\left(e^{\la u,\xi\ra}-1-\la \pi_{J\cup\{k\}}u,\pi_{J\cup\{k\}}h(\xi)\ra\right) M_k(d\xi)\,,\label{eq:R}
\end{eqnarray}
 for $k=1,\ldots,d$ where here we take as truncation function $h(x)=x\mathbbm1_{\{|x|\leq1\}}$. The set of parameters
\begin{equation}\label{eq:admiss}
(b,\beta,a,\alpha,c,\gamma,m,M)\,
\end{equation}
is specified by
\begin{itemize}
 \item $b,\beta_i\in\R^d$ for $i=1,\ldots,d$,
\item $a,\alpha_i\in S^d_{+}$ for $i=1,\ldots,d$, where $S^d_+$ denotes the cones of positive semidefinite $d\times d$ matrices,
\item $c,\gamma_i \in\Rp{}$ for $i=1,\ldots,d$,
\item $m, M_i$ for $i=1,\ldots,d$ are L\'evy measures.
\end{itemize}
\end{theorem}

This set of parameters is called \emph{admissible} if the conditions in Table \ref{tab:admpar} are satisfied with $I$ and $J$ defined as
$I=\{1,\ldots,m\}$ and $J=\{m+1,\ldots,d\}.$ The set of admissible parameters fully characterizes an affine process in $D.$

\begin{table}[!ht]
\begin{displaymath}
 \begin{array}{|ll|}
\hline
\mbox{diffusion}&\\
\hline
        a_{kl}=0 & \mbox{for}\; k\in I\;\mbox{ or}\; l\in I\,,\\
	\alpha_j=0 & \mbox{for all} \;j\in J\,,\\
	(\alpha_i)_{kl}=0 & \mbox{if}\; k\in I\setminus\{i\}\;\mbox{ or }\;l\in J\setminus\{i\}\,,\\
&\\
\hline\hline
\mbox{drift}&\\
\hline
	b\in D,& \\
        (\beta_i)_k\geq0&\mbox{for all}\;i\in I\;\mbox{and}\;k\in I\setminus\{i\}\,,\\
	(\beta_j)_k=0 & \mbox{for all}\; j\in J,\,k\in I\,,\\
&\\
\hline\hline
\mbox{killing}&\\
\hline
	\gamma_j=0&\mbox{for all}\;j\in J,\qquad\qquad\qquad\\
&\\
\hline\hline
\mbox{jumps}&\\
\hline
	\supp m\subseteq D&\mbox{ and }\int_{D\setminus\{0\}}\left((|\pi_I\xi|+|\pi_J\xi|^2)\wedge1\right)m(d\xi)<\infty\,,\\
	M_j=0 & \mbox{  for all }j\in J\,,\\
	\supp M_i\subseteq D& \mbox{ for all } i\in I\mbox{ and}\\
        &\int_{D\setminus\{0\}}\big((|\pi_{I\setminus\{i\}}\xi|+|\pi_{J\cup\{i\}}\xi|^2)\wedge1\big)M_i(d\xi)<\infty\,.\\
\hline
\end{array}
\end{displaymath}
\caption{Set of conditions for admissible parameters.}\label{tab:admpar}
\end{table}

\begin{remark}

If, additionally, the semigroup of transition functions $(p_t)_{t\geq0}$ is homogeneous in the space variable, meaning that, for all $ x\in D\mbox{ and }B\in\Bcal(D)$
$$p_t(x,B)=p_t(0, B-x)\,,$$
then necessarily $R=0$ and it holds
$$\expvBig{x}{e^{\la u,X_t\ra}}=\int e^{\la u,\xi\ra}p_t(x,d\xi)=e^{ t F(u)+\la x,u\ra}\,,$$
for all $(t,x)\in\Rp{}\times D$ and $u\in\Ucal$. Hence $X$ is a (possibly killed) L\'evy process with L\'evy exponent $F$ starting from $x$.
\end{remark}

\subsection{Towards the multivariate Lamperti transform}
When $D=\Rp{}$, it has been proved that there exists a one-to-one correspondence between affine processes
taking values in $D$ and L\'evy processes, see \cite{caballero_lamperti-type_2013}. More precisely, let $Z^{(1)}=(Z^{(1)}_t)_{t\geq0}$ be a L\'evy process starting from 0
taking values in $\R$
whose L\'evy measure has support $\Rp{}$. This implies that there exists a function $R:\im \R\to\C$ such that
$$\expvBig{0}{e^{uZ^{(1)}_s}}=e^{sR(u)}\,,$$
for all $(s,u)\in\Rp{}\times\im\R$. Due to the restrictions on the jump measure, the function $R$ takes the form
$$R(u)=\beta u+\frac{1}{2}\alpha^2 u^2-\gamma+\int_{\Rp{}}\left(e^{u\xi}-1-u\xi\mathbbm1_{\{|\xi|\leq1\}}\right)M(d\xi)\,,$$
where $u\in\im\R$, $\alpha,\beta\in\R$ and $M$ is a measure on $\Rp{}$ which satisfies
$$\int(1\wedge|\xi|^2)M(d\xi)<\infty\,.$$ Moreover, let $Z^{(0)}$ be an independent subordinator with
$$\expvBig{0}{e^{uZ^{(0)}_s}}=e^{sF(u)}\,,$$
for all $(s,u)\in\Rp{}\times\im\R$. Since $Z^{(0)}$ is a subordinator, there exists a constant $b\in\Rp{}$ and a measure $m$ in $\Rp{}$ satisfying
$$\int(1\wedge|\xi|)m(d\xi)<\infty\,,$$
such that, for all $u\in\im \R$,
$$F(u)=b u+\int_{\Rp{}}\left(e^{u\xi}-1\right)m(d\xi)\,.$$
Theorem 2 in  \cite{caballero_lamperti-type_2013} shows that there exists a solution of the following time--change equation
$$X_t=x+Z^{(0)}_t+Z^{(1)}_{\int_0^t X_s ds}$$
for all $(t,x)\in\Rp{}\times \Rp{}$. Moreover, it is proved that the solution is a time homogeneous Markov process, taking values in $\Rp{}$ starting from $x$,
such that the logarithm of the
characteristic function of the transition semigroup is given by an affine function of the initial state $x$.
Hence, by definition, it is an affine process taking values in $\Rp{}$.

Here we are interested in the multivariate generalization of this result, whose weak version is already known in the literature:

\begin{theorem}[Theorem 3.4 in \cite{didactic_2006}]\label{teo:kallsen} Let $X$ be an affine
process with admissible parameter satisfying
$$\int_{\{|\xi|\geq1\}}|\xi_k| M_i(d\xi)<\infty\quad\mbox{and}\quad c=0\,,\gamma_i=0,\quad\mbox{for}\quad 1\leq i, k\leq m.$$
On a possibly enlarged probability space, there exist $d+1$ independent L\'evy processes $Z^{(k)}$ such that
\begin{equation}\label{eq:kallform}
X_t\stackrel{d}{=}x+Z^{(0)}_t+\sum_{k=1}^d Z^{(k)}\left(\int_0^t X^{(k)}_s ds\right)\qquad t\geq0\,.
\end{equation}
\end{theorem}

This result  has to be understood in distributional sense, because, without any additional assumptions, it is not clear how to conclude that the process $X$ is adapted with respect to the (properly time--changed) filtration generated by the L\'evy processes.

In this paper we provide a strong solution of \eqref{eq:kallform} defined on the probability space
$(\Omega,\Gcal,\Pa)$ which carries $Z^{(0)},\ldots,Z^{(d)}$.


\subsection{Affine processes of Heston type}\label{sec:affH}

In this section, we are going to specify a particular subclass of affine processes, which we will call \emph{affine processes of Heston type}.
They are characterized by more restrictive admissible parameters but, at the same time, they constitute a canonical family, in the
sense that every affine process can be obtained as a pathwise transformation of a canonical one. Instead of stating directly the conditions, we work through an example, where we point the main motivations and reasonings for the forthcoming Assumptions \ref{ass:marry} -- \ref{ass:cons}.

\begin{example}
Let us start by writing \eqref{eq:kallform} componentwise. Denote by $Z^{(k,j)}$ the $j$-th coordinate of the $k$-L\'evy process. Then \eqref{eq:kallform} reads
\begin{equation}
\begin{aligned}
X^{(1)}_t&=x_1+Z^{(0,1)}_t+\sum_{k=1}^d Z^{(k,1)}\left(\int_0^t X^{(k)}_s ds\right)\,,\quad t\geq0\\
\vdots&\\
X^{(d)}_t&=x_d+Z^{(0,d)}_t+\sum_{k=1}^d Z^{(k,d)}\left(\int_0^t X^{(k)}_s ds\right)\,,\quad t\geq0\,.\\
\end{aligned}
\end{equation}
Due to the drift conditions summarized in Table \ref{tab:admpar}, we conclude that, for $k=m+1,\ldots,d$, $Z^{(k)}$ is a L\'evy process with triplet $(\beta_k,0,0)$ with $\pi_I\beta_k$ identically zero. Therefore we can write
\begin{equation*}
\begin{aligned}
X^{(1)}_t&=x_1+Z^{(0,1)}_t+\sum_{k=1}^m Z^{(k,1)}\left(\int_0^t X^{(k)}_s ds\right)\,,\quad t\geq0\\
\vdots&\\
X^{(d)}_t&=x_d+Z^{(0,d)}_t+\sum_{k=1}^m Z^{(k,d)}\left(\int_0^t X^{(k)}_s ds\right)+\sum_{k=m+1}^d (\beta_k)_d\left(\int_0^t X^{(k)}_s ds\right)\,,\quad t\geq0\,.\\
\end{aligned}
\end{equation*}

We first transform the process into another affine process with functional characteristic $F=0$. We will see that, up to an enlargement of the state space, there is no loss of generality in assuming that the parameters in the L\'evy--Khintchine form of $F$ are all identically zero. Just for simplicity assume that $n=m=1$. Augment the process $X=(X^{(1)},X^{(2)})$ by considering $$ Y=(Y^{(0)},Y^{(1)},Y^{(2)}):=(1,X^{(1)},X^{(2)})\,.$$ Moreover define, for $k=0,1,2$,
$$\overline Z^{(k)}=(\overline Z^{(k,0)},\overline Z^{(k,1)},\overline Z^{(k,2)}):=(0, Z^{(k,1)},Z^{(k,2)})\,.$$
Then we can write
\begin{displaymath}
  \left(
\begin{array}{c}
Y^{(0)}_t\\
Y^{(1)}_t\\
Y^{(2)}_t
\end{array}
\right)
=  \left(
\begin{array}{c}
1\\
x_1\\
x_2
\end{array}
\right)
+\overline Z^{(0)}\left(\int_0^t Y^{(0)}_s ds\right)+\overline Z^{(1)}\left(\int_0^t Y^{(1)}_s ds\right)+\overline Z^{(2)}\left(\int_0^t Y^{(2)}_s ds\right)\,,\quad t\geq0\, .
\end{displaymath}
Observe that the process $Y$ takes values in $\Rp{2}\times\R$. Hence, up to a change of the state space, we are led to consider solutions of
\begin{equation*}
\begin{aligned}
X^{(1)}_t&=x_1+\sum_{k=1}^m Z^{(k,1)}\left(\int_0^t X^{(k)}_s ds\right)\,,\quad t\geq0\\
\vdots&\\
X^{(d)}_t&=x_d+\sum_{k=1}^m Z^{(k,d)}\left(\int_0^t X^{(k)}_s ds\right)+\sum_{k=m+1}^d (\beta_k)_d\left(\int_0^t X^{(k)}_s ds\right)\,,\quad t\geq0\,.\\
\end{aligned}
\end{equation*}
In order to additionally simplify the system, we introduce a second pathspace transformation, which allows us to work only with affine processes with admissible parameters satisfying the additional property $(\beta_j)_k=0$ for all $j,k\in J$. This means that the L\'evy processes $Z^{(k)}$ with $k=m+1,\ldots,d$ are not only deterministic but actually identically equal to zero. This transformation has been introduced in \cite{keller-ressel_affine_2011} and it is based on the method of the moving frames. The general case will be treated in the proof of Proposition \ref{prop:gocanon}, here we present the case $n=m=1$. Let $X=(X^{(1)},X^{(2)})$ be an affine process in $\Rp{}\times\R$. Consider the process $Y=(Y^{(1)},Y^{(2)})$ with $Y_0=x$ and
\begin{align*}
Y^{(1)}_t&:=X^{(1)}_t-\int_0^t X^{(1)}_s ds\,,\quad t\geq0\\
Y^{(2)}_t&:=X^{(2)}_t-(\beta_2)_2\int_0^t X^{(2)}_s ds\,,\quad t\geq0\,.
\end{align*}
Theorem 5.1 in \cite{keller-ressel_affine_2011} guarantees that $Y$ is again an affine process in $\Rp{}\times\R$ with admissible parameter $\beta^Y_2=(0,0)$. Moreover this transformation can be inverted.  Hence, when $n=m=1$,  up to an invertible pathspace transformation, we can restrict ourselves to the solution of a system of type
\begin{equation*}
\begin{aligned}
X^{(1)}_t&=x_1+ Z^{(1,1)}\left(\int_0^t X^{(1)}_s ds\right)\,,\quad t\geq0\\
X^{(2)}_t&=x_2+ Z^{(1,2)}\left(\int_0^t X^{(1)}_s ds\right)\,,\quad t\geq0\,,\\
\end{aligned}
\end{equation*}

or more generally

\begin{equation*}
\begin{aligned}
X^{(1)}_t&=x_1+\sum_{k=1}^m Z^{(k,1)}\left(\int_0^t X^{(k)}_s ds\right)\,,\quad t\geq0\\
\vdots&\\
X^{(d)}_t&=x_d+\sum_{k=1}^m Z^{(k,d)}\left(\int_0^t X^{(k)}_s ds\right)\,,\quad t\geq0\,.\\
\end{aligned}
\end{equation*}

Then, it is evident that only the equation determining $\pi_I X$ is a real time--change equation. As soon as we provide a strong solution for the system of time--change equations describing the positive components, we automatically find a solution for the components taking values in $\R^n$.
\end{example}

Let $X$ be an affine process taking values in $D$ and denote by $(b,\beta,a,\alpha,c,\gamma,m,M)$ its set of admissible parameters. The next condition implies that the function
$\phi$ in the definition of affine property is identically zero.
 \begin{assumption}\label{ass:marry}
 The condition ${\bf A^\marry}$ is satisfied if $(b,a,c,m)=(0,0,0,0)$.
 \end{assumption}

The next assumption implies that the process is homogeneous in the last $n$ variables.

\begin{assumption}\label{ass:heston}
 The condition $\boldsymbol{A^H}$ is satisfied if, for all $i,j\in J$ it holds $(\beta_i)_j=0$.
 \end{assumption}

Finally, to ensure that a solution of the system exists for all $t\geq0$, we introduce this last set of conditions.
\begin{assumption}\label{ass:cons}
 The condition $\boldsymbol{\mathring A}$ is satisfied if, for all $i\in I$ it holds $c=0\,,\gamma_i=0$ and
 $$\int\left(|\pi_I \xi|\wedge|\pi_I \xi|^2\right)M_i(d\xi),\quad\mbox{for all}\quad i\in I\, .$$
 \end{assumption}

\begin{definition}\label{def:hestonaff}
 We call an affine process with admissible parameters $(b,\beta,a,\alpha,c,\gamma,m,M)$ satisfying ${\bf A^\marry}$, ${\bf A^H}$ and ${\bf \mathring A}$
an \emph{affine process of Heston type}.
\end{definition}

Among the previous assumptions, only Assumption \ref{ass:cons} is a real restriction on the structure on the admissible parameters.  As observed also in \cite{didactic_2006} (also compare with Lemma 9.2. in \cite{2003}) Assumption \ref{ass:cons} guarantees that the solution process does not explode in finite time and hence the time--change process is always well defined. Up to an enlargement of the state space and a pathwise
transformation, there is no loss of generality in assuming that both Assumption \ref{ass:marry} and Assumption \ref{ass:heston} hold.

In the following proposition we present all steps which allow us to reduce a general affine process into an affine process of Heston type.

 \begin{proposition}\label{prop:gocanon} Let $X$ be an affine process satisfying Assumption \ref{ass:cons}.
 On a possibly enlarged probability space, there exists a process $X^\marry$ such that
 \begin{enumerate}
 \item $X^\marry$ is an affine process taking values in $\Rp{m+1}\times\R^n$ satisfying the following property:
  there exists a function $\psi^\marry:\C^{m+1}_{\leq 0}\times\im \R^n\to\C^{d+1}$ such that for all $(t,x^\marry)\in \Rp{}\times (\Rp{m+1}\times\R^n)$ it holds
 $$\expvBig{x^\marry}{e^{\la u,X^\marry\ra}}=e^{\la x^\marry,\psi^\marry(t,u)\ra}$$
 for $u\in\C^{m+1}_{\leq0}\times\im\R^n$,
 \item for all $u=(u_1,u_2)\in \C^{m+1}_{\leq 0}\times\im \R^n$ it holds
 $$\pi_{\{m+2,\ldots,d+1\}}\psi^\marry(t,u_1,u_2)=u_2\,,$$
 for all $t\geq0$,
 \item the set of admissible parameters for $X^\marry$ satisfies the assumptions $\mathbf{A^\marry}$, $\mathbf {A^H}$ and $\mathbf{\mathring A}$ and moreover, for all $k=1,\ldots,m+1$, the matrix
 $\alpha_k$ has the form
\begin{displaymath}
 \alpha_k=\left(
\begin{array}{ccccccc|ccc}
&&&0&&&&&&\\
&&&\vdots &&&&&&\\
&&&0&&&&&&\\
0&\ldots&0&(\alpha_k)_{kk}&0&\ldots&0&0&\ldots&0\\
&&&0&&&&&&\\
&&&\vdots &&&&&&\\
&&&0&&&&&&\\
\hline
&&&0&&&&&&\\
&&&\vdots&&&&&\alpha_k^J&\\
&&&0&&&&&&\\
\end{array}
\right)\,,
\end{displaymath}
with $(\alpha_k)_{kk}\geq0$ and $\alpha_k^J\in S^+_n$,
 \item for all $(t,x)\in\Rp{}\times D$ and $u\in\Ucal$, define $x^\marry=(1,x)$ and $v=(0,u)$. It holds
 $$\expvBig{x}{e^{\la u,X_t\ra}}=\expvBig{x^\marry}{e^{\la v,X^\marry_t\ra}}\,.$$
  \end{enumerate}
  \proof
  Given two indices $i,j=1,\ldots,d+1$ with $i<j$ denote by
$$[i:j]:=\{i,i+1,\ldots,j-1,j\}\,.$$  We start using Proposition 1.23 in \cite{nico_thesis}. Fix $x_0\in\Rp{}$ and define
\begin{align}
& x^\marry:=(x_0,x)\in\Rp{m+1}\times\R^n\,,\\
&\Ucal^\marry:=\C^{m+1}_{\leq 0}\times\im\R^n\,,\\
&\psi^\marry(t,u_0,u_1,\ldots,u_d):=
\left(
 \begin{array}{c}
\phi(t,u_1,\ldots,u_d)+u_0\\
\psi(t,u_1,\ldots,u_d)\\
 \end{array}
\right).\label{eq:Psiinf}
 \end{align}
 Due to regularity in $t$ of $\phi(t,u)$ and $\psi(t,u),$ we conclude that $\psi^\marry(t,\cdot)$ is a regular semiflow. Hence, from Proposition 7.4 in \cite{2003}, we conclude that
there exists an affine process $X^\marry$ with state space $\Rp{m+1}\times \R^n$ satisfying
$$\expvBig{x^\marry}{e^{\la u,X^\marry_t\ra}}=e^{\la x^\marry,\psi^\marry(t,u)\ra}, \quad u\in\Ucal^\marry\, .$$
Now we can apply the method of moving frames (see Theorem 5.1 in \cite{keller-ressel_affine_2011}) to the affine process $X^\marry$. Let
$(0,\beta,0,\alpha,0,0,0,M)$ be its set of admissible parameters. Denote by $\Bcal$ the $d\times d$ matrix obtained by placing each
$\beta_i,\;i=1,\ldots,d$ as a column
\begin{equation}\label{eq:driftmat}
 \Bcal=\left(
\begin{array}{c|c}
 \Bcal_I & 0\\
\hline
\Bcal_{IJ} & \Bcal_J
\end{array}
\right).
\end{equation}
 Define the matrix
\begin{displaymath}
 T=\left(
\begin{array}{c|c}
 I&0\\
\hline0&\Bcal^\top_J
\end{array}
\right)\in\R^{d\times d}
\end{displaymath}
and the map
\begin{eqnarray*}
 \Tcal:X^\marry&\mapsto& X^\marry-T^\top\int_0^\cdot X^\marry_sds\,.
\end{eqnarray*}
The process $\Tcal X^\marry$ is an affine process with Fourier--Laplace transform given by
$$\expvBig{x^\marry}{e^{\la u,\Tcal X^\marry_t\ra}}=e^{\la \pi_{[1:m+1]}x^\marry,\pi_{[1:m+1]}\psi^\marry(t,u)\ra+\la \pi_{[m+2:d+1]} x^\marry,\pi_{[m+2:d+1]}u\ra}.$$
In particular, the assumptions ${\bf A^\marry}$ and ${\bf A^H}$ are satisfied. Now we move on the structure of the matrices $\alpha_k$, $k=1,\ldots,m+1$. Due to the
restrictions on the admissible parameters $\alpha_1$ is already in the specified form  with $(\alpha_1)_{11}=0$. The matrices $\alpha_k$, $k=1,\ldots,m+1$ can be transformed
simultaneously into block diagonal form by means of a linear map. See \cite{filipovic2009affine}.  Finally, if $v=(0,u)$ with $u\in\Ucal$

  $$\expvBig{(1,x)}{e^{\la v,X^\marry_t\ra}}=e^{\phi(t,u)+\la x,\psi(t,u)\ra}=\expvBig{x}{e^{\la u,X_t\ra}}\,,$$
within
\begin{align}
&\psi^\marry(t,(0,u)):=
\left(
 \begin{array}{c}
\phi(t,u)\\
\psi(t,u)\\
 \end{array}
\right)\,.
 \end{align}

  \endproof
 \end{proposition}

\section{Existence of the solution of the time--change equation}\label{sec:existstc}

\subsection{The setting} Let $Z^{(1)},\ldots,Z^{(d)}$ be $d$ independent c\`adl\`ag $\R^d$-valued L\'evy processes, each of them with L\'evy triplet
$(\beta_k,\alpha_k,M_k),\,k=1,\ldots,d$, defined on the same probability space $(\Omega,\Gcal,P)$. Henceforth, we assume that the following restrictions on the L\'evy triplets hold:
\begin{description}
 \item[(H)] the family $(0,\beta,0,\alpha,0,0,0,M)$ consisting of the collection of the triplets $(\beta_k,\alpha_k,M_k)$,\\$k=1,\ldots,d$ satisfies the assumptions
${\bf A^H}$ and ${\bf \mathring{A}}$
\end{description}

Now we consider the process $Z=(Z^{(1)},\ldots,Z^{(d)})\in\R^{d^2}$ on the product space
$$(\Omega,\Gcal,P):=(\prod_{k=1}^d \Omega^{(k)},\otimes_{k=1}^d \Gcal^{(k)},\otimes_{k=1}^d P^{(k)})\, .$$

We fix $x\in D$ and consider the functions
\begin{equation}
f^{(k)}_i(y):=\la x+Ny,e_k\ra,\;\mbox{ for  }k=1,\ldots,d,\;i=1,\ldots,d,\;y\in\R^{d^2}\,,
\end{equation}
where $N\in\R^{d\times d^2}$ is the matrix obtained by horizontally concatenating $d$ times the identity matrix of dimension $d$. In the next section it will be essential to construct the solution of a system of time--change equations of type
\begin{equation}\label{eq:tccoord}
 Y^{(k)}_i(t):=Z^{(k)}_i\left(\int_0^t f^{(k)}_i(Y_s)ds\right),\quad  k,i=1,\ldots,d\,,\mbox{and}\,, t\geq0\,.
\end{equation}

The aim of this section is to prove the following result

\begin{theorem}\label{teo:uniquensol}
Let $Z^{(1)},\ldots,Z^{(d)}$ be $d$ independent $\R^d$-valued L\'evy processes with c\`adl\`ag paths defined on the same probability space $(\Omega,\Gcal,P)$. For  $k=1,\ldots,d$, denote by
$(\beta_k,\alpha_k,M_k)$ the respective L\'evy triplets. Under the assumption that the triplets satisfy {\bf (H)},
for all $x\in D$, there exists a solution of the following time--change problem
\begin{equation}\tag{\ref{eq:tccoord}}
 Y^{(k)}_i(t):=Z^{(k)}_i\left(\int_0^t f^{(k)}_i(Y_s)ds\right),\quad k,i=1,\ldots,d\mbox{ and }t\geq0\,,
\end{equation}
with
\begin{equation}
f^{(k)}_i(y):=\la x+Ny,e_k\ra,\quad k,i=1,\ldots,d,\;y\in\R^{d^2}
\end{equation}
and $N\in\R^{d\times d^2}$ the matrix obtained by horizontally concatenating $d$ times the identity matrix of dimension $d$.
\end{theorem}

\subsection{The proof}

The proof of Theorem \ref{teo:uniquensol} is done in several steps. We first translate the problem of existence and
uniqueness of a solution for \eqref{eq:tccoord} in the problem of existence and uniqueness of a system of ODEs.

Introduce
\begin{equation}
\tau^{(k)}_i(t):=\int_0^t f^{(k)}_i(Y_s)ds,\;\mbox{for }k,i=1\ldots,d\mbox{ and }t\geq0\,,
\end{equation}
 and define
\begin{equation}\label{eq:tauvec}
\underline\tau(t):=(\tau^{(1)}_1(t),\ldots,\tau^{(1)}_d(t),\ldots,\tau^{(k)}_i(t),\ldots,\tau^{(d)}_1(t),\ldots,\tau^{(d)}_d(t))\, .
\end{equation}

Existence of a solution of \eqref{eq:tccoord} is equivalent to the existence of a solution of the following system of ODEs
\begin{equation}\label{eq:ODEtauall}
\left\{
\begin{aligned}
  \dot\tau^{(k)}_i(t)&=f^{(k)}_i\left(\underline Z(\underline\tau(t)\right),\quad\mbox{for all }k,i=1,\ldots,d\,,t\geq0\,,\\
\tau^{(k)}_i(0)&=0\,,
\end{aligned}
\right.
\end{equation}
where
$$\underline Z(\underline\tau(t)):=\left(Z^{(1)}_1(\tau^{(1)}_1(t)),\ldots,Z^{(k)}_i(\tau^{(k)}_i(t)),\ldots, Z^{(d)}_d(\tau^{(d)}_d(t))\right)\,.$$

We start by showing that, existence for a solution of \eqref{eq:ODEtauall} can be proved by focusing on the components with $k=1,\ldots,m$ and $i=1,\ldots,m$.
\begin{lemma}\label{lem:reduction}
 If
\begin{equation}
\left\{
\begin{aligned}
 \dot\tau^{(k)}_i(t)&=f^{(k)}_i(\underline Z(\underline\tau(t))\,,t\geq0\,,\\
\tau^{(k)}_i(0)&=0\,,
\end{aligned}
\right.
\end{equation}
admits a solution for all $k=1,\ldots,m$ and $i=1,\ldots,m$, then it admits also a solution for all $k=1,\ldots,d$ and $i=1,\ldots,d$.
\proof
For all $k=1,\ldots,d$ it holds
$$f^{(k)}_1(y)=x_k+\sum_{h=1}^d y^{(h)}_k=f^{(k)}_2(y)=\ldots=f^{(k)}_d(y),\quad\mbox{for all }y\in\R^{d^2}\,,$$
and therefore
$$\tau^{(k)}_1(t)=\ldots=\tau^{(k)}_d(t),\quad\mbox{for all }t\geq0\,.$$
Denote
$$\tau^{(k)}(t):=\tau^{(k)}_1(t),\mbox{ for all }t\geq0\mbox{ and }k=1,\ldots,d\,.$$
By definition, for each $k=1,\ldots,d\,,$
\begin{align*}
\tau^{(k)}(t)&=\int_0^t \left(x_k+\sum_{i=1}^d Y^{(i)}_k(s)\right)ds\\
&=\int_0^t \left(x_k+\sum_{i=1}^d Z^{(i)}_k(\tau^{(i)}_k(s))\right)ds\\
&=\int_0^t \left(x_k+\sum_{i=1}^d Z^{(i)}_k(\tau^{(i)}(s))\right)ds\,.
\end{align*}
By Assumption \ref{ass:heston}, each $Z^{(j)}$ for $j\in J$ is a L\'evy process with L\'evy triplet $(0,0,0)$, hence it is identically zero.
In particular, once we find a solution for the system
\begin{align*}
\tau^{(k)}(t)&=\int_0^t \left(x_k+\sum_{i=1}^m Z_k^{(i)}(\tau^{(i)}(s))\right)ds\,\quad\mbox{for all }k=1,\ldots,m\,,
\end{align*}
then, for all $t\geq0$ and $j=m+1,\ldots,d$, we can compute
\begin{align*}
\tau^{(j)}(t)&=\int_0^t \left(x_j+\sum_{i=1}^m Z^{(i)}_j(\tau^{(i)}(s))\right)ds\, ,
\end{align*}
since the right hand side of the last equation does not
depend anymore on the left hand side.
\endproof
\end{lemma}

Henceforth, $Z^{(1)},\ldots,Z^{(m)}$ are $m$ independent L\'evy processes on $\R^m$, each of them with L\'evy triplets
$(\beta_i,\alpha_i,M_i)$, $i=1\ldots,m$, satisfying
\begin{displaymath}
 \begin{array}{ll}
	(\alpha_i)_{kl}=0 & \mbox{for all}\; k,l\in I\mbox{ such that } (k,l)\neq(i,i)\,,\\
     (\beta_i)_k\geq0&\mbox{for all}\;i\in I\;\mbox{and}\;k\in I\setminus\{i\}\,.
\end{array}
\end{displaymath}

From Lemma \ref{lem:reduction}, we known that it suffices to study existence and uniqueness of the solution of the following problem
\begin{equation}
\left\{
\begin{aligned}
 \dot\tau^{(k)}(t)&=x_k+\Zcal_k(\underline\tau(t)),\qquad k=1,\ldots,m\,,t\geq0\,,\\
\tau^{(k)}(0)&=0\,,
\end{aligned}
\right.
\end{equation}
where
\begin{eqnarray}
 \Zcal:\Rp{m}&\to&\R^m\nonumber\\
\underline s&\mapsto&\sum_{i=1}^m Z^{(i)}({s_i})\label{def:Zcal}\,.
\end{eqnarray}
In vector notation the previous initial value problem reads
\begin{equation}\label{eq:ODEform}
\left\{
\begin{aligned}
 \dot{\underline\tau}(t)&=x+\Zcal(\underline\tau(t))\,,t\geq0\,,\\
\underline\tau(0)&=0\,.
\end{aligned}
\right.
\end{equation}

\begin{remark}\label{rem:nonegpath}
By definition it holds
$$\dot{\underline\tau}(0)=x\in\Rp{m}\,.$$
Due to the restrictions on the parameters, each $Z^{(i)}$, $i=1,\ldots,m$, is a process with no negative jumps. This implies that, whenever a component $\tau_{i^*}$
reaches zero for some $i^*$, the corresponding
component L\'evy process $Z^{(i^*)}$ is stopped. In particular, each trajectory of $x+\Zcal$ stays positive until it is absorbed at zero.
\end{remark}

\subsubsection{Approximation of the vector field}

In order to construct a solution for \eqref{eq:ODEform}, we seek for a decomposition of type
$$\Zcal=\til{\Zcal}{}+\notil{\Zcal}{}$$
such that the system
\begin{equation}
\left\{
\begin{aligned}
 \dot{\underline\tau}(t)&=(x+\til{\Zcal}{})(\underline\tau(t))\,,t\geq0\\
\underline\tau(0)&=0\,,
\end{aligned}
\right.
\end{equation}
 reduces to a decoupled system of $m$ one dimensional problems and $\notil{\Zcal}{}:=\Zcal-\til{\Zcal}{}$.

The L\'evy--It\^o decomposition, together with the canonical form of the admissible parameters, gives
\begin{align*}
Z^{(i)}_t=&\beta_i t+\sigma_i B^{(i)}_t+\int_0^t \int \xi\mathbbm1_{\{|\xi|>1\}} \Jcal^{{(i)}}(d\xi, ds)\\
&+\int_0^t \int \xi\mathbbm1_{\{|\xi|\leq1\}} (\Jcal^{{(i)}}(d\xi, ds)-M_i(d\xi)ds)
\end{align*}
where $\sigma_i=\sqrt{(\alpha_i)_{ii}}$, $B^{(i)}$ is a process in $\R^m$ which evolves only
along the $i$-th coordinate as Brownian motion and $\Jcal^{(i)}$ is the jump measure of the process $Z^{(i)}$.

Now, from the assumption on the set of admissible parameter,
$$\pi_{I\setminus\{i\}}\beta^{}_i\in\Rp{m-1}\quad\mbox{and}\quad(\beta^{}_i)_i\in\R\,.$$
Decompose

$$Z^{(i)}=:\til{Z}{(i)}+\notil{Z}{(i)}$$
where $\til{Z}{(i)}$ and $\notil{Z}{(i)}$ are two stochastic processes on $\R^m$ defined by
\begin{displaymath}
 \begin{array}{rll}
\til{Z}{(i)}_k(t)&:=0\,,&\mbox{for }k\neq i\,,\\
\til{Z}{(i)}_i(t)&:=\sigma_i B^{(i)}_i(t)+(\beta^{}_i)_it+\int_0^t \int \xi_i\mathbbm1_{|\xi|>1}\Jcal^{{(i)}}(d\xi, ds)\,,&\\
&+\int_0^t \int \xi_i\mathbbm1_{|\xi|\leq 1}(\Jcal^{{(i)}}(d\xi, ds)-M_i(d\xi)ds)\,,&\\
\notil{Z}{(i)}(t)&:=\notil{\beta}{}_it+\int_0^t \int (\xi-\xi_ie_i)\mathbbm1_{\{|\xi|>1\}}\Jcal^{{(i)}}(d\xi, ds)\,,&\\
&+\int_0^t \int (\xi-\xi_ie_i)\mathbbm1_{\{|\xi|\leq1\}}(\Jcal^{{(i)}}(d\xi, ds)-M_i(d\xi)ds)\,,&
 \end{array}
\end{displaymath}
where
$$\notil{\beta}{}_i=\beta_i-e_i(\beta_i)_i\,.$$

The following lemma, which is an obvious consequence of the restrictions on the admissible parameters,
collects some path properties of the processes $\til{Z}{(i)}$ and $\notil{Z}{(i)}$. We would like to remark that both c\`adl\`ag property and this special structure of the paths are essential ingredients of our proof.
\begin{lemma}
For all $i = 1,\ldots,m$ it holds
\begin{enumerate}
 \item $\til{Z}{(i)}$ is a L\'evy process with no negative jumps,
\item  $\notil{Z}{(i)}$ is a process with increasing paths.
\end{enumerate}
\end{lemma}

Introduce, for all $\underline s\in\Rp{m}$,
\begin{align*}
\til{\Zcal}{}(\underline s)&:=\sum_{i=1}^m \til{Z}{(i)}({s_i}),\\
\notil{\Zcal}{}(\underline s)&:=\sum_{i=1}^m \notil{Z}{(i)}({s_i})\,.\\
\end{align*}

We will consider separately the initial value problems with vector fields
$$x+\til{\Zcal}{}\quad\mbox{and}\quad\notil{\Zcal}{}\,.$$

The next result shows that it is possible to find a unique solution for the initial value problem
\begin{equation*}
\left\{
\begin{aligned}
 \dot{\underline\tau}^{}({(t_0,\tau_0,x)};t)&=(x+\til{\Zcal}{})(\underline\tau^{}({(t_0,\tau_0,x)};t)\,,\\
\underline\tau^{}({(t_0,\tau_0,x)};t_0)&=\tau_0\,.
\end{aligned}
\right.
\end{equation*}
Later, we will show how to construct a solution of the general problem.

\begin{proposition}\label{prop:Etaudiag}
There exists a unique solution of
\begin{equation}\label{eq:ODEformlamp}
\left\{
 \begin{aligned}
  \dot{\underline\tau} ({(t_0,\tau_0,x)};t)&=(x+\til{\Zcal}{})(\underline\tau ({(t_0,\tau_0,x)};t)),\\
\underline\tau ({(t_0,\tau_0,x)};t_0)&={\tau_0}\,,
 \end{aligned}
\right.
\end{equation}
with ${\tau_0}\in\Rp{m}$ and $t\geq0$.
\proof
Observe that \eqref{eq:ODEformlamp} is a decoupled system of $m$ equations of type
\begin{equation}\label{eq:ODEformlampd1}
\left\{
 \begin{aligned}
 \dot\tau_i ((t_0,\tau_0,x);t)&=(x_i+\til{Z}{(i)}_i)(\tau_i ((t_0,\tau_0,x);t)),\quad i=1,\ldots,m\,,\\
\tau_i ((t_0,\tau_0,x);t_0)&=\pi_{\{i\}}\tau_0\,.
 \end{aligned}
\right.
\end{equation}
where each $\til{Z}{(i)}_i$ is a L\'evy process with no negative jumps. The existence of a unique solution of \eqref{eq:ODEformlampd1} follows from Section 6.1 in \cite{ethier_markov_1986}.
\endproof
\end{proposition}

For the proof of the general result, we will need to approximate $\notil{\Zcal}{}$ with piecewise constant functions. Fix $M\in\N$ and consider the partition $$\Tcal_M:=\left\{\;\frac{k}{2^{M}},\quad k\geq0\right\}\,.$$
Define the following approximations on the partition $\Tcal_M$:

\begin{align*}
 ^{\uparrow}\notil{Z}{(i,M)}_t&:=\sum_{k=0}^\infty \notil{Z}{(i)}_{k/2^M}\mathbbm1_{[\frac{k}{2^M},\frac{k+1}{2^M})}(t)\,,\\
&\\
 {^{\downarrow}\notil{Z}{(i,M)}_t}&:=\sum_{k=0}^\infty \notil{Z}{(i)}_{(k+1)/2^M}\mathbbm1_{[\frac{k}{2^M},\frac{k+1}{2^M})}(t)\,.\\
\end{align*}
Introduce, for $\underline s\in\Rp{m}$, the processes $ {^{\uparrow}{\notil{\Zcal}{(M)}}}(\underline s)$ and $ ^{\downarrow}{\notil{\Zcal}{(M)}}(\underline s)$ obtained
by taking the sums of $ ^{\uparrow}\notil{Z}{(i,M)}_{s_i}$ and $ ^{\downarrow}\notil{Z}{(i,M)}_{s_i}$ respectively.

\begin{notation}\label{not:augpartition}
 Let
$$\Sigma_M:=\bigcup_{i=1}^m\{s\geq0\st \Delta Z^{(i)}_s>0\}$$
and augment the partition $\Tcal_M$ with $\Sigma_M$. Denote the family obtained in this way by $\Tcal^\Sigma_M$.
\end{notation}

We will first construct a solution for the equation \eqref{eq:ODEform} when $\notil{\Zcal}{}$ is replaced by ${ {^{\uparrow}{\notil{\Zcal}{(M)}}}}$.

Hereafter, given $x,y\in\R^m$, we write $x\leq y\;\mbox{ if }\;x_i\leq y_i,\;\mbox{ for all }\;i=1,\ldots,m\,.$

\subsubsection{The algorithm} Let $\til{\Zcal}{}$ and $ {^{\uparrow}{\notil{\Zcal}{(M)}}}$ be defined as above.
\begin{description}
\item[Input:]
Start by defining the random variables
\begin{align}
\overleftarrow{\sigma}&:=(0,\ldots,0),\\
\overrightarrow{\sigma}(\omega)&:=(\sigma^{(1,M)}_1(\omega),\ldots,\sigma^{(m,M)}_1(\omega)),
\end{align}
where each $\sigma^{(i,M)}_1(\omega)$ is the first jump in the path $t\mapsto { ^{\uparrow}\notil{Z}{(i,M)}_t(\omega)}$.
 \item[Step 1:] Let $\underline\tau((t_0,\tau_0,x);t)$ be the solution of the system \eqref{eq:ODEformlamp} starting from
$$t_0=0,\quad \tau_0=(0,\ldots,0)\quad\mbox{ and }\quad x\in\Rp{m}.$$ Consider the solution of \eqref{eq:ODEformlamp} for all times $t$ such that
\begin{equation}\label{eq:stopcond}\tag{\dag}
\underline\tau((t_0,\tau_0,x);t) < \overrightarrow{\sigma}\,.
\end{equation}
Let $t^*$ be the first time such that \eqref{eq:stopcond} does
not hold anymore. Stop the solution $\underline\tau((t_0,\tau_0,x);\cdot)$ at time $t^*$. Observe that the condition \eqref{eq:stopcond} is violated if there exists an index
$i^*\in\{1,\ldots,m\}$ such that $$\tau_{i^*}((t_0,\tau_0,x);t^*)=\sigma^{(i^*,M)}_1 \, .$$
Notice here that there might be more than one $ i^* $, where the above equality is valid, however, for the sake of convenience,
we assume that there exists only one index for the moment. We will deal with the general case in the proof of Theorem \ref{teo:Etauapprox}.
\item[Step 2:] Update
\begin{align}
\overleftarrow{\sigma}&:=(0,\ldots,\sigma^{(i^*,M)}_1,\ldots,0),\\
\overrightarrow{\sigma}&:=(\sigma^{(1,M)}_1,\ldots,\sigma^{(i^*,M)}_{2},\ldots,\sigma^{(m,M)}_1),\\
x&:=x+\Delta{ ^{\uparrow}\notil{\Zcal}{(M)}}(\overleftarrow{\sigma}),
\end{align}
where $\sigma^{(i^*,M)}_2(\omega)$ is the second jump in the path $t\mapsto{ ^{\uparrow}\notil{Z}{(i^*,M)}_t}(\omega)$.
\item[Step 3:] Let $\underline\tau((t_1,\tau_1,x_1);t)$ be the solution of the system \eqref{eq:ODEformlamp} starting from the \emph{updated} values
$$t_1=t^*,\quad \tau_1=\underline\tau((t_0,\tau_0,x_0);t^*)\quad\mbox{ and }\quad x_1=x\in\Rp{m}\,.$$
As before, we let $\underline\tau((t_1,\tau_1,x_1);\cdot)$ evolve until
\begin{equation}
\underline\tau((t_1,\tau_1,x_1);t^*)< \overrightarrow{\sigma}
\end{equation} holds. As soon as this condition does not holds anymore, we stop again the solution.
\item[End:] Do iteratively \textbf{Step 2} and \textbf{Step 3}.
\end{description}

The above algorithm describes the guiding principle for the proof of the next result:

\begin{theorem}\label{teo:Etauapprox}
There exists a solution of
\begin{equation}\label{eq:ODEformapprox}
\left\{
 \begin{aligned}
  \dot{\underline\tau}^{(M)}({(0,0,x)};t)&=(x+\til{\Zcal}{}+ {^{\uparrow}{\notil{\Zcal}{(M)}}})(\underline\tau^{(M)}({(0,0,x)};t)),\\
\underline\tau^{(M)}({(0,0,x)};0)&={0}\,.
 \end{aligned}
\right.
\end{equation}
\proof
We already did all the main steps for the proof of this result.
Let $\Tcal_M$ and $\Sigma$ be the sets defined in Notation \ref{not:augpartition}. Recall that $\Tcal^\Sigma_M$ is a countable family.
Enumerate the elements in $\Tcal^{\Sigma}_M$ such that
$\sigma^{(i)}_k$ denotes the $k$-th jump of $ ^{\uparrow}\notil{Z}{(i,M)}$. Fix $x\in D$ and set
$$(t_0,\tau_0,x):=(0,0,x)$$
and
\begin{align*}
\overleftarrow{\sigma}&:=(0,\ldots,0),\\
\overrightarrow{\sigma}&:=(\sigma^{(1,M)}_1,\ldots,\sigma^{(i,M)}_{1},\ldots,\sigma^{(m,M)}_1)\,,\\
\end{align*}
where $\sigma^{(i,M)}_{k}$ denotes the $k$-th jump in the path $t\mapsto {^{\uparrow}\notil{Z}{(i,M)}_t}$ for all $i=1,\ldots,m$.
By definition $ {^{\uparrow}{\notil{\Zcal}{(M)}}}(\underline s)=0$ for all $\underline s<\overrightarrow{\sigma}$. Proposition \ref{prop:Etaudiag}
gives the existence of the solution of \eqref{eq:ODEformlamp} with this set of input parameters. Denote it by
$\underline\tau{((t_0,\tau_0,x);t)}$. As soon as the solution $\underline\tau{((t_0,\tau_0,x);t)}$ reaches a jump time for $ {^{\uparrow}{\notil{\Zcal}{(M)}}}$,
the vector field in the equation \eqref{eq:ODEformapprox} changes. Precisely, denote by
$$t_1:=\sup\{t>0\st \underline\tau((t_0,\tau_0,x);t)<\overrightarrow{\sigma}\}\,.$$
Again there might be one or more indices $i^*$, where the condition fails. Collect them in a set $I^*\subseteq\{1,\ldots,m\}$.
Update the values

\begin{equation}\label{eq:updatesigma}
\begin{aligned}
\pi_{I^*}\overleftarrow{\sigma}&:=\pi_{I^*}\overrightarrow{\sigma},\\
\pi_{I^*}\overrightarrow{\sigma}&:=\pi_{I^*}\overrightarrow{\sigma}_{++},\\
\end{aligned}
\end{equation}
where $\overrightarrow{\sigma}_{++}$ contains the next jumps of $ ^{\uparrow}\notil{Z}{(i,M)}$ for all $i\in I^*$ after $\overrightarrow\sigma_i$.
 Then define
\begin{align*}
\tau_1&:=\underline\tau((t_0,\tau_0,x);t_1)\,,\\
x_1&:=x+\Delta{ ^{\uparrow}\notil{\Zcal}{(M)}}(\overleftarrow{\sigma}).
\end{align*}

Now, consider again the solution of \eqref{eq:ODEformlamp}, but this time with parameters $(t_1,\tau_1,x_1)$. Denote it by $\underline\tau{((t_1,\tau_1,x_1);t)}$ and
observe that it is well defined until all the coordinates of $\tau{((t_1,\tau_1,x_1);t)}$ stay below the next jump times of $ {^{\uparrow}{\notil{\Zcal}{(M)}}}$.
We obtain the solution of \eqref{eq:ODEformlamp} by pasting a finite amount of solutions obtained in the time subintervals defined by
$\Tcal^{\Sigma}_M$. Define iteratively, for all $n\geq 1$,
\begin{eqnarray}
 t_{n+1}&:=&\sup\{t>t_n\st\underline\tau((t_n,\tau_n,x_n);t_n)<\overrightarrow{\sigma}\},\\
\tau_{n+1}&:=&\underline\tau((t_n,\tau_n,x_n);t_{n+1}),\\
x_{n+1}&:=&x_n+\Delta{^{\uparrow}\notil{\Zcal}{(M)}}(\overleftarrow{\sigma}),
\end{eqnarray}
where, at each step, $\overleftarrow{\sigma}$ and $\overrightarrow{\sigma}$ are updated using the prescription in \eqref{eq:updatesigma}.
 Continuity follows by construction.
\endproof
\end{theorem}

Now that we have found a solution for the approximated problems, we would like to show convergence to the solution of \eqref{eq:ODEform}.

The following results focus on monotonicity and convergence of \eqref{eq:ODEformapprox}.

\begin{lemma}\label{lem:compareapprox}
Let $i=1,\ldots,m$ and $M\in\N$ be fixed. Then, for all $t\geq0$ it holds
$${ ^{\uparrow}\notil{Z}{(i,M)}_t}\leq \notil{Z}{(i,M)}_t\leq  {^{\downarrow}\notil{Z}{(i,M)}_t}\mbox{ almost surely }.$$
Moreover, for each $\omega\in\Omega$, the sequences $\{{ ^{\uparrow}\notil{Z}{(i,M)}}(\omega)\}_{M\in\N}$ and $\{ ^{\downarrow}\notil{Z}{(i,M)}(\omega)\}_{M\in\N}$ are
monotone in the sense that, for all $t\geq 0$,
$${ ^{\uparrow}\notil{Z}{(i,M+1)}_t}(\omega)\geq { ^{\uparrow}\notil{Z}{(i,M)}_t}(\omega)\,$$
and
$${ ^{\downarrow}\notil{Z}{(i,M+1)}_t}(\omega)\leq {^{\downarrow}\notil{Z}{(i,M)}_t}(\omega)\,.$$

\proof
Since $Z^{(i)}$ has no negative jumps, and, by assumption $(\beta_i)_{k}\geq0$ for all $k\neq i$, the paths of $\notil{Z}{(i,M)}$ are increasing. Therefore,
$$\notil{Z}{(i,M)}_t\geq \notil{Z}{(i,M)}_{k/2^M}= ^{\uparrow}\notil{Z}{(i,M)}_t,\quad\mbox{a.s. for all }t\in[\frac{k}{2^M},\frac{k+1}{2^M})\,.$$
For the same reason,
$$\notil{Z}{(i,M)}_t\leq \notil{Z}{(i,M)}_{{(k+1)}/{2^M}}={ {^{\downarrow}\notil{Z}{(i,M)}_t}},\quad\mbox{a.s. for all }t\in[\frac{k}{2^M},\frac{k+1}{2^M})\,.$$ Now,
since for every $M\in\N$ the
partition $\Tcal_{M+1}$ is obtained by halving all the subintervals in
the partition $\Tcal_M$, it clearly holds
\begin{displaymath}
{ ^{\uparrow}\notil{Z}{(i,M+1)}_t}(\omega)=\left\{
\begin{array}{ll}
{ ^{\uparrow}\notil{Z}{(i,M)}_t}(\omega),&\mbox{ for all }t\in\left[\frac{2k}{2^{M+1}},\frac{2k+1}{2^{M+1}}\right)\,,\\
\\
\notil{Z}{(i,M)}_{(2k+1)/2^{M+1}}(\omega),&\mbox{ for all }t\in\left[\frac{2k+1}{2^{M+1}},\frac{2(k+1)}{2^{M+1}}\right)\,.
\end{array}
\right.
\end{displaymath}
Using again the increasing property of the paths of $\notil{Z}{(i,M)}$ we conclude that
$${ ^{\uparrow}\notil{Z}{(i,M)}_t}\geq \notil{Z}{(i,M)}_t,\quad\mbox{a.s. }$$ because
 $$\notil{Z}{(i,M)}_{(2k+1)/2^{M+1}}\geq \notil{Z}{(i,M)}_{2k/2^{M+1}}=\notil{Z}{(i,M)}_{k/2^{M}}\,. $$
 The case with ${ ^{\downarrow}\notil{Z}{(i,M)}}$ goes analogously.
\endproof
\end{lemma}

\begin{proposition}\label{prop:montaucheck}
Let $M\in\N$ be fixed and denote by ${\underline\tau}^{(M)}((0,0,x);t)$ the solution of \eqref{eq:ODEformapprox} constructed in
Theorem \ref{teo:Etauapprox}. Then, for all $t\geq0$ and $x\in \Rp{m}$ it holds
$${\underline\tau}^{(M)}((0,0,x);t)\leq{\underline\tau}^{(M+1)}((0,{0},x);t),\mbox{ almost surely}\,. $$
\proof This follows by construction using the monotonicity proved in Lemma \ref{lem:compareapprox}. Indeed,
denote by $\Tcal^{\Sigma}_{M}:=\{\sigma^{(M)}_k\}_{k\in\N}$ and $\Tcal^{\Sigma}_{M+1}:=\{\sigma^{(M+1)}_k\}_{k\in\N}$ the set of jump times for
${ ^{\uparrow}{\notil{Z}{(M)}}}$ and ${ ^{\uparrow}{\notil{Z}{(M+1)}}}$ respectively. By construction $\Tcal^{\Sigma}_{M}\subset\Tcal^{\Sigma}_{M+1}$ in the sense that,
for each $\sigma^{(M)}_k\in \Tcal^{\Sigma}_{M}$ there exists $h\in\N$ such that $\sigma^{(M)}_k=\sigma^{(M+1)}_h\in \Tcal^{\Sigma}_{M+1}$.
Denote by $\{\sigma^{(M+1)}_{k_h}\}_{h\in\N}$ the jump times of ${ ^{\uparrow}{\notil{\Zcal}{(M+1)}}}$ occurring on the subinterval
$[\sigma^{(M)}_k,\sigma^{(M)}_{k+1}].$ By construction, there is only one jump inside this interval. Write
$\{\sigma^{(M+1)}_{k_h}\}_{h=1,\ldots,3}$ with $\sigma^{(M+1)}_{k_1}=\sigma^{(M)}_k$ and $\sigma^{(M+1)}_{k_3}=\sigma^{(M)}_{k+1}$
Then $\underline\tau^{(M+1)}$ is obtained by pasting a finite number of solutions of initial value problems with piecewise linear vector field.
For each $h=1,2,3$, ${ ^{\uparrow}{\notil{Z}{(M+1)}}}(\sigma^{(M+1)}_{k_h})\geq{ ^{\uparrow}\notil{Z}{(M)}}(\sigma^{(M)}_{k})$. Therefore,
on each subinterval $[\sigma^{(M)}_k,\sigma^{(M)}_{k+1}]$, the solution
$\underline\tau^{(M+1)}((t_k,\tau_k,x_k);t)$ is constructed by pasting a finite number of solutions of type
$\underline\tau((t_{k,h},\tau_{k,h},x_{k,h});t)$ where $x_{k,h}$ is increasing sequence in $h$.
Hence we conclude that, for all $k\in\N$ and for $t\in[\sigma^{(M)}_k,\sigma^{(M)}_{k+1}]$ it holds
$$\underline\tau^{(M+1)}({(t_k,\tau_k,x_k)};t)\geq\underline\tau^{(M)}({(t_k,\tau_k,x_k)};\sigma^{(M)}_k)\,.$$

\endproof
\end{proposition}

The mast monotonicity argument we need follows directly from the definition of the ODEs
\begin{lemma}
Let $M, t_0, \tau_0$ be fixed and $x\leq y$.
Consider the systems
\begin{equation*}
\left\{
 \begin{aligned}
  \dot{\underline\tau}^{(M)}({(t_0,\tau_0,x)};t)&=(x+\til{\Zcal}{}+ {^{\uparrow}{\notil{\Zcal}{(M)}}})(\underline\tau^{(M)}({(t_0,\tau_0,x)};t)),\\
\underline\tau^{(M)}({(t_0,\tau_0,x)};t_0)&={\tau_0}\,.
 \end{aligned}
\right.
\end{equation*}
\begin{equation*}
\left\{
 \begin{aligned}
  \dot{\underline\tau}^{(M)}({(t_0,\tau_0,y)};t)&=(y+\til{\Zcal}{}+ {^{\uparrow}{\notil{\Zcal}{(M)}}})(\underline\tau^{(M)}({(t_0,\tau_0,y)};t)),\\
\underline\tau^{(M)}({(t_0,\tau_0,y)};t_0)&={\tau_0}\,.
 \end{aligned}
\right.
\end{equation*}
Then, for all $t\geq t_0$ it holds
$${\underline\tau}^{(M)}((t_0,\tau_0,x);t)\leq{\underline\tau}^{(M)}((t_0,\tau_0,y);t),\mbox{ almost surely}\,. $$
\end{lemma}

Finally, due to monotonicity, we know that the sequence $\underline\tau^{(M)}$ admits a limit. With the next result we show that the limit is actually
finite and, by monotone convergence, it coincides with the solution of \eqref{eq:ODEform}.

\begin{proposition}\label{prop:finalapprox}
For all $t\geq0$ and $x\in\Rp{m}$ the sequence $\underline\tau^{(M)}({(0,0,x)};t)$ converges
$$\lim_{M\to\infty}\underline\tau^{(M)}({(0,0,x)};t)=\underline\tau^{(*)}({(0,0,x)};t)\,$$
and the limit can be identified with the solution of \eqref{eq:ODEform}.
\proof
Let $\underline\tau^{(*)}({(0,0,x)};\cdot)$ be the limit of the sequence $\{\underline\tau^{(M)}({(0,0,x)};\cdot)\}_{M\geq0}$. Since the sequence is a monotone sequence, the convergence is actually uniform.
Observe that the same holds for the limit of the sequence of solutions of the system \eqref{eq:ODEformapprox} when ${ ^{\uparrow}{\notil{Z}{(M)}}}$ is replaced by ${ ^{\downarrow}{\notil{Z}{(M)}}}$.
Applying dominate convergence theorem it follows that $\underline\tau^{(*)}({(0,0,x)};t)$ coincides with the solution of \eqref{eq:ODEform}.
\endproof
\end{proposition}

At this point, most of the results we need for the proof of Theorem \ref{teo:uniquensol} have been proved. The final step is to construct the solution of the
time change equation \eqref{eq:tccoord}
using the solution of the system \eqref{eq:ODEform}.

\paragraph*{Proof of Theorem \ref{teo:uniquensol}}
Let $\underline\tau=(\tau^{(1)},\ldots,\tau^{(m)})$ be the solution of \eqref{eq:ODEform} with $\Zcal:=\sum_{i=1}^m \pi_I Z^{(i)}$. Then,
for $k=1,\ldots,m$ and $i=1,\ldots,d$
$$Y^{(k)}_i(t):=Z^{(k)}_i(\tau^{(k)}(t))\,.$$
Moreover observe that, due to the restrictions in {\bf (H)}, the L\'evy processes $Z^{(k)}$ for $k=m+1,\ldots,d$ are identically zero and therefore also $Y^{(k)}$
are identically zero.
\endproof

\section{Pathwise construction of affine processes with time--change}\label{sec:pathconstr}

We start summarizing the results from Chapter \ref{sec:existstc}. In Proposition \ref{prop:finalapprox} we have shown that the system of ODEs
\begin{equation}\tag{\ref{eq:ODEform}}
\left\{
\begin{aligned}
 \dot{\underline\tau}(t)&=x+\Zcal(\underline\tau(t))\,,\\
\underline\tau(0)&=0
\end{aligned}
\right.
\end{equation}
admits a solution which can be constructed as the limit of approximated problems. Then, in Theorem \ref{teo:uniquensol} we showed how to use these
solutions in order to construct the processes $\{Y^{(k)}_i\}_{k,i=1,\ldots,d}$ defined by means of the time--change equation \eqref{eq:tccoord}.
In this section we are going to see how to combine these processes in order to construct an affine process. Before to do it, we clarify the main
steps by means of an easy two dimensional example with $n=m=1$.
\begin{example}
The results in Section \ref{sec:existstc} in the particular case when $m=n=1$ give the existence of a solution for the time change equation
 $$Y^{(k)}_i(t):=Z^{(k)}_i\left(\int_0^t f^{(k)}_i(Y_s)ds\right),\qquad\mbox{for }k,i=1,2\mbox{ and }t\geq0.$$
Under the assumption {\bf(H)}, $Y^{(2)}_t=0$ for all $t\geq0$ and

$$Y^{(1)}(t):=Z^{(1)}\left(\int_0^t \left(x_1+Y^{(1)}_1(s)\right)ds\right)\,.$$
Define
$$X=x+NY\quad\mbox{with}\quad
N=\left(
\begin{array}{cccc}
1&0&1&0\\
0&1&0&1\end{array}
\right)\,.$$

Inserting the definitions of the $Y^{(k)}_i$, it is clear that $X$ satisfies

\begin{displaymath}
\left(
\begin{array}{c}
 X^{(1)}\\
X^{(2)}
\end{array}
\right)
=\left(
\begin{array}{c}
 x_1\\
x_2
\end{array}
\right)
+
\left(
\begin{array}{c}
 Z^{(1)}_1(\int_0^\cdot X^{(1)}_s ds)+ Z^{(2)}_1(\int_0^\cdot X^{(2)}_s ds)\\
 Z^{(1)}_2(\int_0^\cdot X^{(1)}_s ds)+Z^{(2)}_1(\int_0^\cdot X^{(2)}_s ds)
\end{array}
\right)\,.
\end{displaymath}
In vector notation, we can write
$$X=x+\sum_{i=1}^2  Z^{(i)}\left(\int_0^\cdot X^{(i)}_s ds\right)\,$$
which is indeed the formulation in Theorem \ref{teo:kallsen}.
\end{example}
The next theorem is a re-formulation of the above argument in the general multivariate case. The additional problem we still need to address is measurability of the time--change process with respect to the filtration generated by the $Z^{(i)}$, $i=1,\ldots,d$. In order to do it, we will need the notion of multivariate filtration and multivariate stopping time taken from \cite{ethier_markov_1986}. Recall that $Z^{(i)}$, $i=1,\ldots,d$ is a family of L\'evy processes as in the setting of Section \ref{sec:existstc}.

For all $\underline s=(s_1,\ldots,s_{d^2})\in\Rp{d^2}$ ,
define the $\sigma$-algebra
\begin{equation}\label{eq:multinatural}
\Gcal^\natural_{\underline s}:=\sigma\left(\{Z^{(h)}_{t_h},\; t_h\leq s_h,\;\mbox{for }h=1,\ldots,d^2\}\right)\,,
\end{equation}
and then complete it by
\begin{equation}\label{eq:EK2.7}
\Gcal_{\underline s}=\bigcap_{n\in \N}\Gcal^\natural_{\underline s^{(n)}}\vee\sigma(\Ncal),
\end{equation}
where $\Ncal$ is the collection of sets in $\Gcal$ with $P$-probability zero and $\underline s^{(n)}$ is the sequence defined by $s^{(n)}_k=s_k+1/n$.

\begin{definition}
A random variable $\underline\tau=(\tau_1,\ldots,\tau_{s^2})\in\Rp{d^2}$ is a $(\Gcal_{\underline s})$-stopping time if
$$\{\underline\tau\leq\underline s\}:=\{\tau_1\leq s_1,\ldots,\tau_{d^2}\leq s_{d^2}\}\in\Gcal_{\underline s},\mbox{ for all }\underline s\in\Rp{d^2}\, .$$
\end{definition}

If $\underline \tau$ is a stopping time,
$$\Gcal_{\underline\tau}:=\{B\in\Gcal\st B\cap\{\underline \tau\leq \underline s\}\in\Gcal_{\underline s}\mbox{ for all }\underline s\in\Rp{d^2}\}\, .$$

Now that we have introduced the necessary notation, we are ready to prove the following result.

\begin{theorem}\label{teo:multilam}
Let $(b,\beta,a,\alpha,c,\gamma,m,M)$ be a set of admissible parameters satisfying the Assumptions
${\bf A^\marry}$, ${\bf A^H}$ and ${\bf \mathring A}$.
\begin{enumerate}[(i)]
\item The time--change equation
\begin{equation}\label{eq:timechange}
 X_t=x+\sum_{i=1}^d Z^{(i)}(\theta^{(i)}_t),\;\mbox{ with }\;\theta^{(i)}_{t}=\int_0^t X^{(i)}_{r}dr\, ,
\end{equation}
admits a unique solution.
\item[(ii)]Define
$$\underline \theta^x_{t}:=(\underbrace{\theta^{(1)}_t,\ldots,\theta^{(1)}_t}_{d\,\mbox{\scriptsize times}},\ldots,\underbrace{\theta^{(d)}_t,\ldots,\theta^{(d)}_t}_{d\,\mbox{\scriptsize times}})\in\R^{d^2}\, .$$
The random variable $\underline \theta^x_{t}$ is a $\Gcal_{\underline s}$ stopping time for all $t\geq0$. Hence the time--change filtration
 $$\Gcal_{\underline \theta^x_{t}}:=\{A\st A\cap\{\underline \theta^x_{t}\leq \underline s\}\in\Gcal_{\underline s},\;\mbox{ for all }\underline s\in\Rp{d^2}\}\, ,$$
is well defined.
\item[(iii)] Let $R$ be the function defined as in \eqref{eq:R}.
The solution of \eqref{eq:timechange} is an affine process with functional characteristics $(0,R)$ with respect to the time--changed filtration
$(\Gcal_{\underline\theta^x_{t}})_{t\geq0}$.
\end{enumerate}

\proof

Let $Y\in\R^{d^2}$ be the process obtained by casting the solutions of \eqref{eq:tccoord} as
$$Y:=(Y^{(1)}_1,\ldots,Y^{(1)}_d,Y^{(2)}_1,\ldots,Y^{(2)}_d,\ldots,Y^{(d)}_1,\ldots,Y^{(d)}_d)\, .$$
Consider the matrix
\begin{displaymath}
 N:=\left(\
\begin{array}{ccccccccccc}
 1&0&\cdots&0&1&0&\cdots&0&1&\cdots&0\\
0&1&0&\cdots&0&1&0&\cdots&0&\cdots&0\\
&&\ddots&&&&\ddots&&&\\
0&0&\cdots&1&0&0&\cdots&1&\cdots&0&1\\
\end{array}
\right)\in{\R^{d\times d^2}}\,.
\end{displaymath}
Then
$$X=x+\sum_{k=1}^d Y^{(k)}$$
is a solution of time-change equation \eqref{eq:timechange}. Indeed, in vector notation, we can write
$$X=x+NY\, .$$
Then, if $Z^{(k)}_j$ denotes the $j$-th coordinate of the $k$-th L\'evy process,
\begin{eqnarray*}
 Z^{(k)}_j\left(\int_0^t f^{(k)}_j(Y_s)ds\right)&=&Z^{(k)}_j\left(\int_0^t \la x+NY_s, e_k\ra ds\right)=Z^{(k)}_j\left(\int_0^t X^{(k)}_s ds\right)
\end{eqnarray*}
and
\begin{eqnarray*}
 X_j&=&x_j+\sum_{k=1}^d Y^{(k)}_j\\
&=&x_j+\sum_{k=1}^d Z^{(k)}_j\left(\int_0^t X^{(k)}_s ds\right).
\end{eqnarray*}

Now we move on the measurability of the time--change process. Observe that Theorem \ref{teo:uniquensol} implies that the vector
$$\underline\tau(t):=(\tau^{(1)}_1(t),\ldots,\tau^{(1)}_d(t),\ldots,\tau^{(d)}_1(t),\ldots,\tau^{(d)}_d(t))$$
where
$$\tau^{(k)}_i(t)=\int_0^t f^{(k)}_i(Y_s)ds$$ is a $\Gcal_{\underline s}$ stopping time for all $t\geq0$. This follows from Theorem VI.2.2. in \cite{ethier_markov_1986}. From the affine relationship between $X$ and $Y$ we conclude that
$\underline\theta_t$ is a $\Gcal_{\underline s}$ stopping time and therefore the time--changed filtration is well defined.

Now, we need to check that $X$ is a homogeneous Markov process with respect to $(\Gcal_{\underline\theta^x_{t}})_{t\geq0}$.
Applying Proposition I.6 in \cite{bertoin_levy_1998} at each component $Z^{(k)},\;k=1,\ldots,d$ we get that
$(\underline Z(\underline \theta^x_{t+h})-\underline Z(\underline \theta^x_{t}))_{h\geq0}$
has the same law as $\underline Z(\underline \theta^x_{h})_{h\geq0}$ and it is independent of $\Gcal_{\underline \theta^x_{t}}\,.$

Therefore
$$X^x_{t+h}=X^x_{t}+N\big(\underline Z(\underline\theta^x_{t+h})-\underline Z(\underline\theta^x_t)\big)=:\Scal_t(\underline Z(\underline \theta^x_{t+h})-\underline Z(\underline \theta^x_{t}),X^x_t)\, ,$$
with
\begin{eqnarray*}
\Scal_t:(\R^{d^2},\prod_{i=1}^d(\Gcal_{\underline \theta^x_t}))\times(\R^d,\Gcal_{\underline \theta^x_t})&\to&(\R^d,\Gcal_{\underline \theta^x_t})\\
(Z,X)&\to& X+NZ.
\end{eqnarray*}
Therefore, we conclude that the conditional law of $X^x_{t+h}$, given $\Gcal_{\underline\theta^x_{t}}$, is
$X^x_t$ measurable.
Markov property translates into
$$X^x_{t+h}=\Scal_0(\underline Z(\underline \theta^y_{{h}}),y)_{| y=X^x_t}\, .$$

Additionally the time--change process is absolutely continuous with $$\frac{d}{dt}\theta_i^{(k)}(t)=X^{(k)}_{t^-},\;\mbox{ for all }k,i=1,\ldots,d\, .$$
The characteristics of the time--changed semimartingale can be computed using the formulas in Theorem 8.4. in \cite{barndorff-nielsen_change_2010} from where we conclude that the process
$(\Scal_0(\underline Z(\underline\theta^{x}_{t}),x))_{t\geq0}$
has characteristics $(\beta(X_{-}),\alpha(X_{-}),M(X_{-}))$, where
\begin{align*}
 \beta(x)&=x_1\beta_1+\ldots+x_m\beta_m\,,\\
\alpha(x)&=x_1\alpha_1+\ldots+x_m\alpha_m\,,\\
M(x,B)&=x_1M_1(B)+\ldots+x_m M_m(B),\quad B\in\Bcal(D)\,.
\end{align*}

\endproof
\end{theorem}




\bibliographystyle{alpha}
\bibliography{sample}

\begin{thebibliography}{CPGUB13}

\bibitem[Bat96]{bates_jumps_1996}
D.~S. Bates.
\newblock Jumps and stochastic volatility: exchange rate processes implicit in
  deutsche mark options.
\newblock {\em Review of Financial Studies}, 9(1):69--107, 1996.

\bibitem[Ber98]{bertoin_levy_1998}
J.~Bertoin.
\newblock {\em {L}\'evy {P}rocesses}.
\newblock Cambridge University Press, 1998.

\bibitem[BNS01]{barndorff-nielsen_shepard}
O.~E. Barndorff-Nielsen and N.~Shephard.
\newblock Modelling by {L}\'evy processes for financial econometrics.
\newblock In {\em L\'evy processes}, pages 283--318. Birkh\"auser Boston,
  Boston, MA, 2001.

\bibitem[BNS10]{barndorff-nielsen_change_2010}
O.~E. Barndorff-Nielsen and A.~Shiryaev.
\newblock {\em {C}hange of {T}ime and {C}hange of {M}easure}.
\newblock World Scientific, 2010.

\bibitem[CPGUB13]{caballero_lamperti-type_2013}
M.~E. {C}aballero, J.~L. {P}{\'e}rez {G}armendia, and G.~{U}ribe {B}ravo.
\newblock A {L}amperti-type representation of continuous-state branching
  processes with immigration.
\newblock {\em The Annals of Probability}, 41(3):1585--1627, May 2013.

\bibitem[CT13]{cuchiero_path_2011}
C.~Cuchiero and J.~Teichmann.
\newblock Path properties and regularity of affine processes on general state
  spaces.
\newblock In {\em S\'eminaire de {P}robabilit\'es {XLV}}, volume 2078 of {\em
  Lecture Notes in Math.}, pages 201--244. Springer, Cham, 2013.

\bibitem[DFS03]{2003}
D.~Duffie, D.~Filipovi{\'c}, and W.~Schachermayer.
\newblock Affine processes and applications in finance.
\newblock {\em The Annals of Applied Probability}, 13(3):984--1053, 2003.

\bibitem[EK86]{ethier_markov_1986}
S.~N. Ethier and T.~G. Kurtz.
\newblock {\em Markov processes}.
\newblock Wiley Series in Probability and Mathematical Statistics: Probability
  and Mathematical Statistics. John Wiley \& Sons, Inc., New York, 1986.

\bibitem[FM09]{filipovic2009affine}
D.~Filipovic and E.~Mayerhofer.
\newblock Affine diffusion processes: {T}heory and applications.
\newblock {\em Advanced Financial Modelling}, 8:1--40, 2009.

\bibitem[Gab14]{nico_thesis}
N.~Gabrielli.
\newblock {\em Affine processes from the perspective of path space valued
  L\'evy processes}.
\newblock PhD thesis, ETH Z\"urich, 2014.

\bibitem[Hes93]{heston_closed-form_1993}
S.~L. Heston.
\newblock A closed-form solution for options with stochastic volatility with
  applications to bond and currency options.
\newblock {\em The Review of Financial Studies}, 6(2):327--343, 1993.

\bibitem[Kal06]{didactic_2006}
J.~Kallsen.
\newblock A didactic note on affine stochastic volatility models.
\newblock In {\em From Stochastic Calculus to Mathematical Finance}, page 343.
  Bachelier Colloquium on Stochastic Calculus and Probability, 2006.

\bibitem[KST11]{keller-ressel_affine_2011}
M.~{Keller-Ressel}, W.~Schachermayer, and J.~Teichmann.
\newblock Affine processes are regular.
\newblock {\em Probab. Theory Related Fields}, 151(3-4):591--611, 2011.

\bibitem[KST13]{keller-ressel_regularity_2011}
M.~{Keller-Ressel}, W.~Schachermayer, and J.~Teichmann.
\newblock Regularity of affine processes on general state spaces.
\newblock {\em Electron. J. Probab.}, 18:no. 43, 17, 2013.

\end{thebibliography}







\end{document}